%%%%%%%%%%%%%%%%%%%%%%%%%%%%%%%%%%%%%%%%%%%%%%%%%%%%%%%%%%%%%%%%%%%%%%%%%%
%%%%%%%%%%%%%%%%%%%%%%%%%%%%%%%%%%%%%%%%%%%%%%%%%%%%%%%%%%%%%%%%%%%%%%%%%%
%%%%%%%%%%%%%%%%%%%%%%%%%%%%%%%%%%%%%%%%%%%%%%%%%%%%%%%%%%%%%%%%%%%%%%%%%%
%%%%%%%%%%%%%%%%%%%%%%%%%%%%%%%%%%%%%%%%%%%%%%%%%%%%%%%%%%%%%%%%%%%%%%%%%%
%%%%%%%%%%%%%%%%%%%%%%%%%%%%%%%%%%%%%%%%%%%%%%%%%%%%%%%%%%%%%%%%%%%%%%%%%%
%This is closed_revised003.tex. 
%Preliminaries

\input amstex

\documentstyle{amsppt}

\loadbold

\magnification=\magstep1

\pageheight{9.0truein}
\pagewidth{6.5truein}

%%%%%%%%%%%%%%%%%%%%%%%%%%%%%%%%%%%%%%%%%%%%%%%%%%%%%%%%%%%%%%%%%%%%%%%%%%
%%%%%%%%%%%%%%%%%%%%%%%%%%%%%%%--MACROS--%%%%%%%%%%%%%%%%%%%%%%%%%%%%%%%%%
%%%%%%%%%%%%%%%%%%%%%%%%%%%%%%%%%%%%%%%%%%%%%%%%%%%%%%%%%%%%%%%%%%%%%%%%%%
\def\L{\Lambda}
\def\Mod{\operatorname{Mod}}
\def\End{\operatorname{End}}
%%%%%%%%%%%%%%%%%%%%%%%%%%%%%%%%%%%%%%%%%%%%%%%%%%%%%%%%%%%%%%%%%%%%%%%%%%
%%%%%%%%%%%%%%%%%%%%%--REFERENCE-MACROS--%%%%%%%%%%%%%%%%%%%%%%%%%%%%%%%%%
%%%%%%%%%%%%%%%%%%%%%%%%%%%%%%%%%%%%%%%%%%%%%%%%%%%%%%%%%%%%%%%%%%%%%%%%%%
\def\Pap{1}
\def\Pop{2}
\def\Ros{3}
\def\Smione{4}
\def\Smitwo{5}
\def\StaVdB{6}
\def\Ste{7}
\def\VdB{8}
%%%%%%%%%%%%%%%%%%%%%%%%%%%%%%%%%%%%%%%%%%%%%%%%%%%%%%%%%%%%%%%%%%%%%%%%%%
%%%%%%%%%%%%%%%%%%%%%%--END-OF-MACROS--%%%%%%%%%%%%%%%%%%%%%%%%%%%%%%%%%%%
%%%%%%%%%%%%%%%%%%%%%%%%%%%%%%%%%%%%%%%%%%%%%%%%%%%%%%%%%%%%%%%%%%%%%%%%%%

\topmatter

\title A Remark on Closed Noncommutative Subspaces \endtitle

\rightheadtext{Closed Noncommutative Subspaces}

\author E. S. Letzter \endauthor

\abstract Given an abelian category with arbitrary products, arbitrary
coproducts, and a generator, we show that the closed subspaces (in the sense
of A. L. Rosenberg) are parameterized by a suitably defined poset of ideals in
the generator. In particular, the collection of closed subspaces is itself a
small poset. \endabstract

\address Department of Mathematics, Temple University, Philadelphia, PA 19122
\endaddress

\email letzter\@math.temple.edu \endemail

\thanks The author thanks the Department of Mathematics at the
University of Pennsylvania for its hospitality; the research
for this paper was undertaken while he was a visitor on sabbatical
there. The author is grateful for support during this period from a
Temple University Research and Study Leave Grant. This research was also
supported in part by a grant from the National Security Agency.
\endthanks

\endtopmatter

\document

\head 1. Introduction \endhead

Following Rosenberg \cite{\Ros}, Van den Bergh \cite{\VdB}, and others
(the reader is referred to \cite{\StaVdB} for an overview), the
fundamental objects of study in noncommutative algebraic geometry are
Grothendieck categories, interpreted as categories of sheaves on (not
explicitly defined) noncommutative schemes. Key to this approach is
the notion of a {\sl closed subspace}, investigated in detail, for
example, in \cite{\Pap; \Ros; \Smione; \Smitwo; \VdB}.  This brief
note records an elementary -- but apparently previously unnoticed --
observation about closed subspaces, working in the slightly more
general setting of an abelian category equipped with a generator and
having arbitrary products and coproducts: We observe that there is a
duality between the collection of closed subspaces and a suitably
defined small poset of {\sl ideals\/} within the generator. (After the
fact, the poset structure does not depend on the choice of generator.)
This duality generalizes Rosenberg's duality, in the case of module
categories, between closed subspaces and two-sided ideals \cite{\Ros,
III.6.4.1, p\. 127}. A different generalization can be found in
\cite{\VdB, \S 3.4}.

As a corollary to the duality we present below, it follows that the collection
of closed subspaces is a small set. Also, the closed subspaces satisfy the
descending chain condition if and only if the poset of ideals satisfies the
ascending chain condition. Consequences of the descending chain condition for
closed subspaces can be found, for example, in \cite{\Pap; \Smione}.

\head 2. Closed Subspaces and the Poset of Ideals \endhead

Our approach is heavily influenced by \cite{\Ros; \Smione; \VdB}. The
reader is referred to \cite{\Pop; \Ste} for background information on
abelian categories.

\subhead 2.0 Notation and Assumptions \endsubhead (i) Throughout this note,
$A$ will denote an abelian category with coproducts and products (both over
arbitrary sets) and equipped with a generator $\L$. Examples for $A$ include
module categories over rings (always assumed to be associative and unital) and
Grothendieck categories. To avoid set-theoretic difficulties, we will always
view $A$ as large and sets as small.

(ii) Recall that subobjects, quotient objects, and subquotient objects (i.e.,
quotients of subobjects) refer to equivalence classes of objects (cf\., e.g.,
\cite{\Ste, pp\. 83--84}). We will abuse the notation slightly by referring to
morphisms, coproducts, products, etc\. of these equivalences classes (rather
than of their representative objects).

\subhead 2.1 \endsubhead Let $M$ be an object in $A$. Since $A$ has a
generator, the collection $L(M)$ of subobjects of $M$ is a set
(cf\. \cite{\Ste, IV.6.6}) and is also a lattice under sums, intersections,
and inclusions (as defined in \cite{\Ste, p\. 88}). In particular, the
collection of subquotients of $M$ is a set.

\subhead 2.2 \endsubhead Following \cite{\Smione, \S 2; \VdB, \S 3.3} we
will say that a full subcategory $C$ of $A$ is {\sl closed\/} if
$C$ is closed under products, coproducts, subquotients, and
isomorphisms. In the terminology of \cite{\Smione}, the closed
subcategories are {\sl closed subspaces}.

\subhead 2.3 \endsubhead Let $C$ be a closed subcategory of $A$. The
collection of kernels of all of the morphisms from $\L$ into objects in $C$ is
a set, and we will use $I_{\L}(C)$, the {\sl ideal of $C$ in $\L$}, to denote
the intersection of these kernels. Also, we will say that a subobject of $\L$
is an {\sl ideal in $\L$\/} if it is the ideal of some closed subcategory of
$A$. 

As of now we do not know of an equivalent, more intrinsic, general description
of the ideals in $\L$.

\subhead 2.4 \endsubhead For rings, the preceding definition coincides with
the standard theory: Let $R$ be a ring, and let $I$ be a (two-sided) ideal of
$R$, defined in the usual ring-theoretic way. Let $\Mod R$ denote the category
of left $R$-modules, with generator $R$, and let $C(I)$ denote the full
subcategory of $\Mod R$ consisting of those left $R$-modules $M$ for which
$I.M = 0$. It is easy to see that $C(I)$ is a closed subcategory of $\Mod
R$. It is also easy to see, in the notation of (2.3), that $I = I_R(C(I))$. In
\cite{\Ros, III.6.4.1} A. L. Rosenberg proved that every closed subcategory of
$\Mod R$ has the form $C(I)$ for some ideal $I$, and it follows that the
ideals defined in (2.3) for the generator $R$ of $\Mod R$ are exactly the
usual ring-theoretic ideals of $R$. Hence, there is a bijective
duality between closed subcategories of $\Mod R$ and ideals of $R$. It is this
duality that we seek to generalize below. A different generalization may be
found in \cite{\VdB, \S 3.4}.

\proclaim{2.5 Proposition} Let $C$ be a closed subcategory of $A$, and set $I
= I_\Lambda(C)$. Then $\L/I$ is contained in $C$ and is a generator for
$C$. Consequently, if $C$ and $C'$ are distinct closed subcategories of $A$,
then $I_\L(C) \ne I_\L(C')$. \endproclaim

\demo{Proof} Let $S$ be the set of all quotients of $\L$ contained in $C$. Let
$M \in C$ be the product of these quotients, and observe that $I$ is equal to
the intersection of the kernels of the canonical quotient morphisms from $\L$
into objects contained in $S$. Therefore, $I$ is the kernel of the resulting
product morphism from $\L$ into $M$, and so $\L/I \in C$. Now note that every
object in $C$ is an epimorphic image of a coproduct of isomorphic copies of
$\L$, since $\L$ is a generator for $A$ (cf\. \cite{\Pop, 2.8.2,
p\. 51}). In other words, every object in $C$ is a sum of images of $\L$. But
the set of images of $\L$ in $C$ coincides with the set of images of $\L/I$ in
$C$, and so every module in $C$ is a sum of images of $\L/I$.  Consequently,
every module in $C$ is an epimorphic image of a coproduct of ismorphic copies
of $\L/I$. Hence, $\L/I$ is a generator for $C$ (again cf\. \cite{\Pop, 2.8.2,
p\. 51}), and the proposition follows. \qed\enddemo

\proclaim{2.6 Corollary} The collection of closed subcategories of $A$ is a
small set. \qed \endproclaim

\subhead 2.7 \endsubhead When $I$
is an ideal of $\L$ we will use $C_\Lambda(I)$ to denote the unique closed
subcategory of $A$ whose ideal in $\L$ is $I$. We regard the set of
ideals in $\L$ as a poset with respect to inclusion.  We also regard the set
of closed subcategories of $A$ as poset with respect to inclusion (i.e.,
$C \subseteq C'$ when every object of $C$ is contained in $C'$.)

\proclaim{2.8 Lemma} Let $C$ and $C'$ be closed subcategories of
$A$ with respective ideals $I$ and $I'$ of $\Lambda$. 
Then $C \subseteq C'$ if and only if $I' \subseteq I$. \endproclaim

\demo{Proof} It follows immediately from the definition of ideal that
$I' \subseteq I$ if $C \subseteq C'$. Now suppose that $I' \subseteq
I$. Then $\L/I$ is a quotient of $\L/I'$, and in particular, $\L/I$ is
contained in $C'$. Since $\L/I$ is a generator for $C$, it follows
that $C \subseteq C'$. \qed\enddemo

We can combine (2.5) and (2.8) to obtain the following:

\proclaim{2.9 Theorem} The function 
%%%%%%%%%%%%%%%%%%%%%%%%%%%%%%%%%%%%%%%%%%%%%%%%%%%%%%%%%%%%%%%%%%%%%%%%%%%%
$$\{ \text{Closed Subcategories of $A$} \} @> \quad C \;
\longmapsto \; I_\L(C) \quad >> \{ \text{Ideals in $\L$} \}$$
%%%%%%%%%%%%%%%%%%%%%%%%%%%%%%%%%%%%%%%%%%%%%%%%%%%%%%%%%%%%%%%%%%%%%%%%%%%%
is a poset isomorphism. \qed\endproclaim

\subhead 2.10 \endsubhead In view of (2.9), we see, after the fact, that the
poset of ideals in $\L$, up to a poset isomorphism, does not depend on the
choice of generator. As an application of this fact, we see that $A$ satisfies
the descending chain condition for closed subcategories if and only if the
poset of ideals satisfies the ascending chain condition. (Consequences of the
descending chain condition for closed subcategories can be found, for example,
in \cite{\Pap; \Smione}.) As another application, when $\Lambda$ is
projective, it then follows from (2.4) and (e.g.) \cite{\Pop, 7.4} that the
poset of ideals of $\L$ is isomorphic to the poset of ring-theoretic ideals in
the ring $\End(\L,\L)$. (We leave for the future a more concrete description
of the ideals in a projective generator.)

\subhead Acknowledgements \endsubhead I am grateful to Peter Jorgensen and
Paul Smith for useful remarks on an earlier draft of this note. I am also
grateful for the suggestions of the referee, which helped to more clearly
present the issues raised in this note.

\enddocument